\documentclass[12pt]{article}
\usepackage{amsmath}
\usepackage{amsthm}
\usepackage{amssymb}
\def\R {{\mathbb R}}
\def\ds{\displaystyle}
\theoremstyle {plain}
\newtheorem{thm}{Theorem}[section]

\newtheorem{prop}[thm]{Proposition}
\newtheorem{cor}[thm]{Corollary}
\theoremstyle{definition}

\theoremstyle{remark}
\newtheorem{rem}{Remark}[section]
\begin{document}

\centerline{\Large\bf Realization of the annihilation operator}

\bigskip

\centerline{\Large\bf for generalized oscillator-like system}

\bigskip

\centerline{\Large\bf by a differential operator
\footnote{This research was supported by RFFI grant No 00-01-00500}}

\bigskip
\bigskip

\centerline{\large\bf Borzov V.V.\quad and Damaskinsky E.V.}

\medskip

\begin{abstract}
This work continues the research of generalized Heisenberg
algebras connected with several orthogonal polynomial systems. The
realization of the annihilation operator of the algebra
corresponding to a polynomial system by a differential operator
$A$ is obtained. The important special case of orthogonal
polynomial systems, for which the matrix of the operator $A$ in
$l^2(Z_+)$ has only off-diagonal elements
 on the first upper  diagonal different from zero, is considered.
The known generalized Hermite polynomials give us an example of
such orthonormal system. The replacement of the usual derivative
by q-derivative allows us to use the suggested approach for
similar investigation of various "deformed" polynomials.
\end{abstract}

\bigskip

\tableofcontents

\section{Introduction}

The connection of the  orthogonal polynomials theory with the
group theory as well as the connection the theory of quantum
(deformed) groups and algebras with the theory of basic
hypergeometric functions (deformed
polynomials)~\cite{Andr,AskWil,Chih,GasRa} is well
known~\cite{Vil,KlVil,KS,kor}. In particular, the Hermite
polynomials system (after multiplication on $\exp(-x^2)$) is a
system of eigenfunctions for  quantum-mechanical harmonic
oscillator~\cite{FDL} energy operator. Besides, these polynomials
are connected with representations of the algebra generated by
Heisenberg commutation relations (and with representations of the
appropriate group as well). Similarly, many known $q$-deformed
Hermite polynomials~\cite{AskWil,asi,flo} make up  eigenfunctions
systems of the energy operators for the related deformed
oscillators~\cite{bie,mac,dam}. Moreover, these polynomials
naturally arise from the analysis of irreducible representations
of appropriate deformed algebras.

At the standard approach it is suggested that a representation of
considered algebra (or group) is prescribed. Then we look for a
system of orthogonal polynomials which make up the eigenfunctions
system of a significant operator of this algebra (for example,
Hamiltonian). In other words, we look for the basis in the
representation space, in which this operator take the diagonal
form.

In \cite{B} we considered in some sense an inverse task . Let
$\mu$ be a positive  Borel measure on the real line $R^1$. We
consider a set of  real orthonormal polynomials in the space
$L^2(R^1;\mu(dx))$. For the orthogonal polynomial system  one
wants to construct oscillator like algebra, such that these
polynomials make up the eigenfunctions system of the energy
operator corresponding to this oscillator. By the full solution of
this problem it is suggested that a differential or a difference
equation for the given system of polynomials is known. In
work~\cite{B} (strictly) classical polynomials, (namely, the
Hermite, Jacobi and Laguerre polynomials) were considered as
examples. From these examples it follows that the central problem
with obtaining the differential (or difference) equation for the
given system of polynomials was in finding  a realization of the
annihilation (lowering) operator by a differential operator. Note
that the operator of a finite order in work~\cite{B}  occurres
only for the standard quantum mechanical oscillator connecting
with the orthonormal Hermite polynomial system. In all other cases
considered in~\cite{B}  the annihilation operator connected with
the orthonormal classical polynomial system was realized by a
differential operator of the infinite order. More exactly these
operators were described as rather simple functions of the first
order differential operator. Other examples of realization of the
annihilation operator by the known functions of the first order
differential operator are considered in works~\cite{AAS,ARS,AFW}.

In the present work we shall obtain a expression for an
annihilation operator of the generalized deformed oscillator
algebra related to a complete system of the real orthonormal
polynomials. The above-mentioned  examples from \cite{B} are
concrete realizations of the general expressions obtained here.

Let us remark that the operators of the finite order arise for the
 Hermite-Chihara polynomial systems
(see, \cite{sze} \cite{Chih} and \cite{Bor})
 and the Hermite-Hahn polynomial systems \cite{H}.

It is more interesting to apply the obtained results to
investigating a generalized oscillator connected with symmetric
basic number $[a]=\ds\frac{q^a-q^{-a}}{q-q^{-1}}$ of quantum
groups theory. This case is still not clearly understood. The
research of this case is in a stage of completion now.

For the convenience of the reader  in the following section we
recall some necessary information from the theory of Jacobi
matrices and the classical  moment problem.

\section{Background information}

Let a operator $X$ is given by action on the standard orthonormal
basis $\left\{e_n | n\in Z_{+} \right\} $ in the Hilbert space
$H=l^2(Z _{+})$:
\begin{equation}\label{b1}
Xe_n =b_n\, e_{n+1} + a_n\, e_n + b_{n-1}\, e_{n-1},
\qquad a_n\geq 0,\ b_n\in{\R}.
\end{equation}
Then $X$ can be represent in $H$ by an infinite 3-diagonal matrix
\begin{equation}\label{b2}
X=\left(
\begin{array}{cccccc}
a_0 & b_0 &  0  &  0  & \cdots  & \cdots  \\
b_0 & a_1 & b_1 &  0  & \cdots  & \cdots  \\
 0  & b_1 & a_2 & b_2 & \cdots  & \cdots  \\
 0  & 0   & b_2 & a_3 & \ddots  & \cdots  \\
\vdots  & \vdots  & \vdots  & \ddots  & \ddots  & \ddots  \\
\vdots  & \vdots  & \vdots  & \vdots  & \ddots  & \ddots
\end{array}
\right) \,,
\end{equation}
known as the Jacobi matrix.

With the Jacobi matrix it is related the set of polynomials
$P_n(x)$ of degree $N$, satisfying the recurrent relation,
\begin{equation}\label{b3}
b_n\, P_{n+1}(x)+a_n\, P_n(x)+b_{n-1}P_{n-1}(x) = xP_n(x)\, ,
\end{equation}
and subjecting to the following "initial conditions"
\begin{equation}\label{b4}
 P_0(x)=1\, ,\qquad P_{-1}(x)=0\, .
\end{equation}

Further, we shall assume that $a_n,\ b_n \in {\R}$ and $a_n=0\ ,
b_n>0$. Under such conditions polynomials $P_n(x)$ have real
coefficients and fulfill the following parity conditions $P_n(-x)=
(-1)^n P_n(x).$
 Notice that the recurrent relations
(\ref{b3}) have two linearly independent solutions. The
polynomials $P_n(x)$, which solve to (\ref {b3}) and subject to
the initial conditions (\ref {b4}), are called  polynomials of the
first kind. Besides, there are independent set of solutions
consisting from polynomials $Q_n(x)$, which satisfy with the
following initial conditions
\begin{equation}
\label{b5} Q_0(x) =0,\qquad\qquad Q_1(x)=\frac{1}{b_0}.
\end{equation}
For these polynomials the standard agreement ${\rm deg}\,Q_n(x)=n-1$
is carried out.  These polynomials are known as polynomials of the second
kind for the Jacobi matrix $X$ (\ref {b2}).

The polynomials of the first and second kind are related by the
following equation
\begin{equation} \label{b6}
P_{n-1}(x)Q_{n}(x) - P_{n}(x)Q_{n-1}(x)=\frac{1}{b_{n-1}},
\qquad (n=1,2,3,...)\, .
\end{equation}

It is known that polynomials
$\left\{ P_n(x)\right\}_{n=0}^{\infty}$
are orthonormal
\begin{equation}\label{b7}
\int_{\R} P_n(x)P_m(x) \, d\sigma(x) = \delta_{n,m}
\end{equation}
with respect to a positive Borel measure $\sigma$ on $\mathbb{R}$.
In order to evaluate the measure $\sigma$ one must solve the
following task.

By  $s_n$ we denote the factor $\alpha_0$ at $P_0(t)$ in the
expansion
$$ t^n=\sum_{k=0}^n\alpha_k\ P_k(t)\ . $$
For the given
numerical sequence $\left\{ s_n \right\} _{n=0}^{\infty}$ it is
required to find a positive measure $\sigma$ on the real line such
that the following equalities:
\begin{equation}\label{b8}
s_n =\int_{-\infty }^\infty t^n d\sigma(t)\, ,
\quad n=0,1,2,\ldots\
\end{equation}
are hold. This task is known as the degree Hamburger moment
problem~\cite{L10}.

If the system (\ref{b8}) defines the measure uniquely, then the
moment problem is called {\bf determined} one. Otherwise there
exists an infinite set of such measures, and the moment  problem
is called {\bf undetermined} one.

The matrix (\ref{b2}) determines a symmetrical operator $X$
in $H.$ It is known that the deficiency indices for the
operator $X$ are equal  $(0,0)$ or $(1,1).$

\begin{prop}
Let $\overline{X}$ be a closure of an operator $X$ in $H$. The
following conditions are equivalent:
\begin {enumerate}
\item the deficiency indices of an operator $X$
are equal $(0,0)$;
\item an operator $\overline{X}$ is selfadjoint in the Hilbert space $H$;
\item the moment problem (\ref {b8}) for
the given number set $\left\{ s_n \right\} _{n=0}^{\infty}$ such
that
\begin{equation}\label{b9}
s_n = \left( e_0, X^n e_0 \right)\,
\end{equation}
is  determined;
\item a series
$\sum_{n=0}^{\infty} \left| P_n(z) \right|^2$
is divergent for all $z\in \mathbb{C}$ such that
${\rm Im}z\neq 0 $.
\end {enumerate}
\end {prop}

\bigskip

Let the operator $X$ has the deficiency indices $(1,1)$. Then the
operator $\overline{X}$ not selfadjoint and there exists infinite
number of its selfadjoint extensions. In this case series
$\sum_{n=0}^{\infty} \left| P_n(z) \right|^2 <\infty$ is divergent
at all $z\in \mathbb{C}$ such that ${\rm Im}z\neq 0 $, and moment
problem  (\ref{b8}) for a set of numbers $s_n $ from (\ref {b9})
is indetermined.

The following proposition gives a sufficient condition under which
the operator $X$ has the deficiency indices $(1,1)$.

\begin {prop}
Under the following conditions
\begin {enumerate}
\item there exists $N$ such that for all $n>N$ the inequality
\begin{equation}\label{b10}
b_{n-1}\ b_{n+1} \leq b_{n}^{\,\,2};
\end{equation}
\item
\begin{equation}\label{b11}
\sum_{n=0}^{\infty}\frac{1}{b_n} < \infty ;
\end{equation}
\end {enumerate}
the operator $X$ has the deficiency indices equal to $(1,1).$
\end{prop}

\bigskip

Moreover, we have the following statement.

\begin {prop} (see \cite {L11}).
If the conditions (\ref {b10}) and (\ref {b11}) are fullfiled,
there is an infinite number of selfadjoint extensions of the operator
$\overline{\mathit {X}}$, and related moment problem (\ref {b8})
for a set of numbers $\left\{ s_n \right\} _{n=0}^{\infty}$ from
(\ref {b9}) is undetermined.
\end {prop}

\bigskip

The polynomials $P_n$ and $Q_n$ can are represented in the
form~\cite{dam}
\begin{gather}
P_n(x;q)=\sum_{m=0}^{\epsilon (\frac n2)}
\frac{(-1)^m}{\sqrt{\left\{ n\right\}!}} b_0^{2m-n}
\alpha _{2m-1,n-1}x^{n-2m}  \label{b12}
\\
\alpha _{\,-1;n-1}\equiv 1;\qquad \alpha
_{2m-1;n-1}=\sum\limits_{k_1=2m-1}^{n-1}\left\{ k_1\right\}
\sum\limits_{k_2=2m-3}^{k_1-2}\left\{ k_2\right\} \ldots
\sum\limits_{k_m=1}^{k_{m-1}-2}\left\{ k_m\right\} ,\quad m\geq 1
\label{b13}
\\
Q_{n+1}(x;q)=\sum_{m=0}^{\epsilon (\frac n2)}
\frac{(-1)^m}{\sqrt{\left\{n+1\right\} !}}
 b_0^{2m-n}\beta _{2m,n}x^{n-2m}  \label{b14}
\\
\beta _{0;n}\equiv 1;\qquad \beta _{2m;n}=
\sum\limits_{k_1=2m-1}^n\left\{k_1\right\}
\sum\limits_{k_2=2m-2}^{k_1-2}\left\{ k_2\right\}
\ldots \sum\limits_{k_m=2}^{k_{m-1}-2}
\left\{ k_m\right\} ,\quad m\geq 1 , \label{b15}
\end{gather}
where $\{s\}=\ds\frac{b_{s-1}^{2}}{b_0^{2}}$, and the integral
part of $x$ is denoted by $\epsilon (x)=\mathrm{Ent}(x)$. Because
the representation  (\ref{b12}) - (\ref{b15}) of the polynomials
$P_{n}$ and $Q_{n}$ have not been adequately explored, we shall
recall some properties of coefficients
$\alpha_{m;n}$ and  $\beta_{m;n}.$

From the recurrent relation (\ref{b3}) it follows that the
coefficients $\alpha_{m;n}$ and $\beta_{m;n}$ satisfy the
conditions:
\begin{align}
\alpha _{2m-1;n}&=\left\{ n\right\}
\alpha _{2m-3;n-2}+\alpha _{2m-1;n-1};
\label{b16} \\
\beta _{2m;n}&=\left\{ n\right\}
\beta _{2m-2;n-2}+\beta _{2m;n-1}.
\label{b17}
\end{align}

From the definitions (\ref{b13}) and (\ref{b15})  of coefficients
$\alpha_{k,n}$ and $\beta_{k,n}$ it follows that
\begin{equation}\label{b18}
\alpha_{2k-1,2k-1}=\{2k-1\}!!,\qquad
\beta_{2k-2,2k-2}=\{2k-2\}!!.
\end{equation}
Here
\[
\left\{ 2n\right\} !!=\left\{ 2n\right\} \left\{ 2n-2\right\}
\cdot \ldots \cdot \left\{ 2\right\} ,\qquad
\left\{ 2n-1\right\} !!=\left\{ 2n-1\right\}
\left\{2n-3\right\} \cdot \ldots \cdot \left\{ 1\right\} .
\]

\section {Statement of a problem}

Let $\frak{H}=L^2\left(\Bbb{R};\mu(\rm{d}x)\right)$ be a Hilbert
space, and $\left\{ \Psi_n(x)\right\}_{n=0}^{\infty}$ be a system
of polynomials, which are orthonormal with respect to the measure
$\mu$. The recurrent relations of these polynomials take the
following  form:
\begin{equation}\label{c1}
b_n\, \Psi_{n+1}(x)+b_{n-1}\Psi_{n-1}(x) = x\Psi_n(x)\, ,\qquad b_n>0\,,
\end{equation}

\begin{equation}\label{c2}
 \Psi_0(x)=1\, ,\qquad \Psi_{-1}(x)=0\, .
\end{equation}

In the work \cite{B} it was shown, that one can construct the
oscillator-like algebra ${\mathcal{A}}_{\Psi}$ corresponding to
this polynomial system. In addition, the space ${\frak H}$ is a
space of the Fock representation, and the polynomials $\left\{
\Psi_n(x)\right\}_{n=0}^{\infty}$ make up the Fock basis. The
generators ${a_\mu}^+,\,{a_\mu}^-,\,N$ of this algebra acts in the
Fock representation as follows
\begin{align}
{a_\mu}^+\Psi_n(x)&=\sqrt{2}b_n{\Psi_{n+1}}(x), &n\geq 0;
\nonumber\\
{a_\mu}^-\Psi_n(x)&=\sqrt{2}b_{n-1}{\Psi_{n-1}}(x),
&n\geq 1\;\quad {a_\mu}^{-}{\Psi_0}(x)=0;\label{c3}\\
N\Psi_n(x)&=n\Psi_n(x), &n\geq 0; \nonumber\\
\end{align}
From these relations it follows that
\begin{equation}
\label{c4} {a_\mu}^-{a_\mu}^+\Psi_n(x)=2{b_n}^{2}\Psi_n(x),\qquad
{a_\mu}^+{a_\mu}^-\Psi_n(x)=
2{b_{n-1}}{}^{2}\Psi_n(x),\qquad n\geq 0,
\end{equation}
and $b_{-1} = 0.$ Then we have the following commutation
relations:
\begin{equation}\label{bd1.11}
\left[ {a_\mu}^-,{a_\mu}^+\right]=2\left(
B(N+1)-B(N)\right),\qquad \left[ N,{a_\mu}^{\pm}\right]=\pm
{a_\mu}^{\pm},
\end{equation}
between the generators of the algebra $\mathcal{A}_{\Psi}.$ The
function $B(N)$ is defined by the formula
\begin{equation}\label{bd1.12}
B(N)\Psi_n(x)=b_{n-1}{}^{2}\Psi_{n-1}(x).
\end{equation}
If the sequence
 $\left\{ b_n\right\}_{n=0}^{\infty}$
satisfies the recurrent relation
\begin{equation}\label{bd1.13}
b_{n}{}^{2}-Qb_{n-1}{}^{2}=C(n),
\end{equation}
then, as follows from (\ref{bd1.12}), the relation
\begin{equation}\label{bd1.13a}
{a_\mu}^{-}{a_\mu}^{+}-Q{a_\mu}^{+}{a_\mu}^{-}=2C(N).
\end{equation}
is fulfilled also.

According to (\ref{b12}) and (\ref{b13}), we have
\begin{equation}
\label{b20}
\Psi_n(x)=\sum_{m=0}^{\epsilon(n/2)} \frac{(-1)^m}{\sqrt{\{n\}!}}
\alpha_{2m-1,\,n-1}b_0^{2m-n}x^{n-2m}, \qquad n\geq0.
\end{equation}

For simplification of  further formulas we  define the operator
$A$ by a relation $a^{-}=\sqrt{2}A.$ Then
\begin{gather}
\label{bd2.3}A\Psi _n=b_{n-1}\Psi _{n-1},\quad
n\geq 1,\qquad A\Psi _0=0;\\
A\Psi _0=0\;\Rightarrow a_{n0}=0,\quad n\geq 0 \label{bd2.4}
\end{gather}

The main purpose of the present work is to find, using the
formulas   (\ref {bd2.3}) and (\ref {bd2.4}), the conditions on
the coefficients $a_{ks}$, under which the operator $A$ take the
form
\begin{equation}
A=\sum_{k,s=0}^\infty a_{ks}x^k\frac{{\rm d^s}}{{\rm d}x^s}.
\label{bd2.2}
\end{equation}

 Substituting (\ref{bd2.2}) and (\ref{b3}) into (\ref{bd2.3}),
one obtains
\begin{multline}
\quad\sum_{k,s=0}^\infty a_{ks}\sum_{m=0}^{\epsilon (n/2)}
\frac{(-1)^m}{\sqrt{\left\{ n\right\} !}} \alpha
_{2m-1,n-1}b_0^{\,2m-n}x^{n+k-2m-s} \frac{\left( n-2m\right)
!}{\left( n-2m-s\right) !}=
\\
b_{n-1}\sum_{m=0}^{\epsilon (\frac{n-1}2)}
\frac{(-1)^m}{\sqrt{\left\{ n-1\right\} !}}
\alpha_{2m-1,n-2}b_0^{\,2m+1-n}x^{n-1-2m}\quad
\label{bd2.5}
\end{multline}
 Equating the coefficients at $x^t$ in both sides of a relation
(\ref{bd2.5}), we receive $$a_{n,0}=0, \qquad  n\geq 0$$
\begin{multline}
\quad\sum_{s=0}^\infty \sum_{m=\delta ( \frac{n-s-1}2)}^{\epsilon
(\frac{n-s}2)}(-1)^ma_{t+s+2m-n,s} \alpha_{2m-1,n-1}b_0^{\,2m-n}
\frac{\left( n-2m\right) !}{\left( n-2m-s\right) !}=
\\
\\
\left\{
\begin{array}{ccc}
0 & if & n-1-t=2p+1,\quad p\geq 0 \\
\sqrt{\left\{ n\right\} }b_{n-1}(-1)^p\alpha _{2p-1,n-2}b_0^{\,-t} & if &
n-1-t=2p,\quad p\geq 0
\end{array}
\right.\,,\quad\quad
\label{bd2.6}
\end{multline}
where the least integer greater than $a$ or equal to $a$ is
denoted by $\delta(a)$.

 From (\ref{bd2.6}) it follows that the problem of a realization
of the annihilation operator by a differential operator
(\ref{bd2.2}) was rather cumbersome in general case. For this
reason we slightly change the statement of a problem. Namely, we
will seek for not an annihilation operator but a "lowering
operator" $A$ (denoted by the same symbol). In other words from
(\ref {bd2.6}) we should look for a operator $A$ such that:
\begin{equation}
\label{bd2.7}A\Psi _n=\gamma_{n}\Psi _{n-1},\quad n\geq 1,\qquad
A\Psi _0=0,
\end{equation}
under the supposition, that the operator (\ref{bd2.2}) works in
space ${H}$ as a shift operator. Thus we supposed that among
coefficients $\left\{ a_{k,s}\right\}_{k,s=0}^{\infty}$ in the
matrix of this operator only elements on the first upper diagonal
are different from zero
\begin{equation}
a_{n,k}=\left\{
\begin{array}{ccc}
0&if&k\neq n+1\\ a_{n,k}\neq 0&if&k=n+1
\end{array}
\right.
\label{bd2.8}
\end{equation}
Let us remark  that condition
\begin{equation}
\label{bd2.9}
a_{n,0}=0,\qquad n\geq 0
\end{equation}
and relations (\ref{bd2.5}) and (\ref{bd2.6})
(with a replacement
 $b_{n-1}\longrightarrow\gamma_n$  in a right hand side)
are hold.

In the following section we shall obtain some consequence of the
relations (\ref{bd2.7}) and (\ref{bd2.8}).

\section {Conditions under which the matrix of an operator $A$ has
non zero elements only upper above diagonal}

Taking into account that all elements $a_{k, s}$ except for
$a_{k,k + 1}$ are equal to zero, and using the relation
\ref{bd2.6} for $t=0,1,2,\ldots$ and $P=0,1,2,\ldots$, we obtain
the following conditions on  coefficients
$\left\{b_n\right\}_{n=0}^{\infty}$
 and
$\left\{\gamma_n\right\}_{n=1}^{\infty}:$
\begin{gather}
\label{bd2.11} \sqrt{\{2p+1\}}\gamma_{2p+1} =
\frac{\alpha_{2p-1,2p}}{\alpha_{2p-1,2p-1}}\cdot\gamma_1,
\qquad p\geq 1\\
\label{bd2.13}\sqrt{\left\{ 2p+2\right\} }\gamma _{2p+2}=
\frac{\alpha_{2p-1,2p+1}}{\alpha _{2p-1,2p}}
\sqrt{\left\{ 2\right\} }\gamma _2,\qquad
p\geq 0
\end{gather}
\begin{equation}
\label{bd2.16} \frac{a_{0,1}}{k!}+ \frac{a_{1,2}}{(k-1)!}+\ldots +
\frac{a_{k-1,k}}{1!}+a_{k,k+1}= \frac{b_0}{(k+1)!}\sqrt{\left\{
k+1\right\} }\gamma _{k+1}
\end{equation}
where $k\geq 0.$ Besides, we have the following consistency
conditions (for $p\geq 0,$ $ k\geq 0$)
\begin{equation}
\label{bd2.17}
\frac{\alpha _{2p+1,2p+k+1}}{\alpha _{2p+1,2p+k}}=
\frac{\alpha_{2p-1,2p+k+1}}{\alpha _{2p-1,2p+k}}\cdot
\frac{\alpha _{1,k+1}}{\alpha _{1,k}}
\end{equation}

As a result we obtain the following theorem.
\begin{thm}
For the operator $A$ defined by (\ref{bd2.2}) fullfil the
conditions (\ref{bd2.7}) and (\ref{bd2.8}) it is necessary and
sufficient to have the conditions (\ref{bd2.11}) and
(\ref{bd2.13}). Here  $\alpha _{2m-1,n-1}$  are the coefficients
of the function
$\Psi_n(x)$ defined by (\ref {bd2.17}). Under these conditions the
elements  $a_{k,k+1}$ of the matrix of the operator $A$ are
determined by the relations (\ref {bd2.16}).
\end {thm}

\begin {cor}
Under the additional conditions
$$
a_{0,1}\neq 0,\qquad 0=a_{1,2}=a_{2,3}=\ldots =a_{k,k+1}=\ldots,
$$
we have
\begin {equation}\label{bd2.18a}
\sqrt{\left\{ k\right\} }\gamma _k=k\gamma_1,\qquad k\geq 1.
\end {equation}
\end {cor}

\begin {thm}
The coefficients  $\alpha_{2m-1,n-1}$ defined by (\ref{b13})
satisfy the conditions (\ref{bd2.17}) iff
there exist a positive sequence
$\left\{ v_k\right\} _1^{\infty},$
such that
\begin{eqnarray}
\alpha _{2m-1,n-1}= \frac{\left\{ 2m-1\right\} !!
\left( v_1v_2\ldots v_{n-1}\right)}
{\left( v_1v_2\ldots v_{2m-1}\right)
\left( v_1v_2\ldots v_{n-2m-1}\right)},
\qquad m\geq 1,\;n\geq3\nonumber\\
\alpha _{-1,n-1}=1,\qquad \alpha_{1,1}=\left\{ 1\right\} =1.
\label{bd2.18}
\end{eqnarray}
 Under these conditions the sequence
$\left\{ v_k\right\} _1^{\infty}$
has the following properties
\begin {enumerate}
\item $1=v_0\leq v_1 \leq v_2\leq\ldots$
\item
$v_{n-2}v_{2p-1}+v_{2p-3}v_{n-2p}=v_{n}v_{2p-3}+v_{2p-1}v_{n-2p}$
for all
$n\geq 2,$ $p\geq 1,$ $2p\leq n$ $(v_{-1}\equiv 0)$.
\end {enumerate}
\end {thm}
\begin {proofname}.

The sufficiency of the condition (\ref{bd2.18}) is checked by the
direct substitution (\ref{bd2.18}) into (\ref{bd2.17}).

Necessity. Let the coefficients $\alpha _{2m-1,n-1}$ determined by
relations (\ref{b13}) are satisfied the conditions (\ref{bd2.17}).
Let us put $v_0=a_{1,1}=1,\quad v_1\geq{1}$, and define
the remaining terms of a sequence $\left\{ v_k\right\}_1^{\infty}$
from the equalities
\begin{equation}\label{bd2.19}
\alpha _{1,n}=\frac{v_{n-1}v_n}{v_1},\qquad n\geq{1}.
\end{equation}
Using  relations (\ref{bd2.17}), we immediately  check the
validity of the formulas (\ref{bd2.18}).
\end{proofname}

An example of a polynomial system of the above-mentioned type give
us the generalized Hermite polynomials  (entered in one of notes
in the monography Сеге~\cite {sze} and explicitly investigated
in~\cite {Chih}). We consider these polynomials in the next
section.

\section {Generalized Hermite polynomials}

The more complete exposition of the material considered in this
section can be found in \cite{Bor}.

Let us consider the Hilbert space
\begin{equation}
{\tt H}_{\gamma}=
{L^2}(R;|x|^{\gamma}(\Gamma(\frac{1}{2}
(\gamma+1)))^{-1}\exp(-x^2){dx}),
\qquad\gamma\geq{-1}.
\label{triv40}
\end{equation}
According to the method suggested in \cite{B}, we shall construct
a canonical orthonormal  polynomial system
$\left\{{\psi_{n}(x)}\right\}_{n=0}^{\infty }$, which is complete
 in the space ${\tt H}_{\gamma}$. The polynomials
$\psi_{n}(x)$ will satisfy the recurrent relations (\ref{c1}) and
(\ref{c2}) with coefficients
$\left\{{b_{n}}\right\}_{n=0}^{\infty }$,
which are given by the following formulas:
\begin{equation}
b_{n-1}=\frac{1}{2} \Biggl\{
\begin{array}{cc}
\sqrt{n}&{n=2m},\\[.3cm]
\sqrt{n+\gamma}&{n=2m+1}.\\
\end{array}
\Biggr. \label{triv42}
\end{equation}
The polynomials  $\psi_{n}(x)$ are defined by the relations
(\ref{bd2.18}),  (\ref{b18}) and the formulas (\ref{b12}), where
 $\{s\}=\ds\frac{b_{s-1}^{2}}{b_0^{2}}$. Besides, the sequence
 $\left\{ {v_{n}}\right\}_{n=0}^{\infty }$
is given by the formulas
\begin{equation}
v_n= \Biggl\{
\begin{array}{cc}
\ds\frac{\gamma+n+1}{\gamma+1}&{n=2m},\\[.5cm]
\ds\frac{n+1}{\gamma+1}&{n=2m+1}.\\
\end{array}
\Biggr. \label{triv41}
\end{equation}
It is clear that
$v_0=1,\quad{v_1=\ds{\frac{2}{\gamma+1}}=b_0^{-2}}.$

The appropriate "lowering" operator $A$ is defined by relations
(\ref{bd2.2}) and (\ref{bd2.8}), where
\begin{equation}
a_{0,1}=1, \qquad
a_{m-1,m}={\frac{(-2)^{m-1}}{m!}}{\frac{\gamma}{\gamma+1}} ,
\qquad m\geq2 . \label{triv43}
\end{equation}
The sequence $\left\{ {\gamma_{n}}\right\}_{n=0}^{\infty}$
included in (\ref{bd2.7}) is given by  the expressions:
\begin{equation}
\gamma_{n}=\frac{\sqrt{2}}{\gamma+1} \Biggl\{
\begin{array}{cc}
\sqrt{n}&{n=2m},\\ \sqrt{n+\gamma}&{n=2m+1}.\\
\end{array}
\Biggr. \label{triv44}
\end{equation}
The annihilation operator $a_{\mu}^{-}$ of the generalized
oscillator algebra is equal to
\begin{equation}
a_{\mu}^{-}=\frac{\gamma+1}{\sqrt{2}}A. \label{triv45}
\end{equation}

 It is valid the following commutation relations for the
position operator $X_\mu$ and the number operator $N$:

\begin{equation}
X_{\mu}\frac{d}{dx}-N=(a_{\mu}^{-})^{2}. \label{triv49}
\end{equation}
Using the eigenvalue equation for a Hamiltonian
$H_{\mu}=a_{\mu}^{-}a_{\mu}^{+}+a_{\mu}^{+}a_{\mu}^{-}$,
and taking into account the equality
$a_{\mu}^{+}a_{\mu}^{-}=2B(N)$,
one can obtain the following relation :
\begin{equation}
a_{\mu}^{-}a_{\mu}^{+}=2B(N+I). \label{triv52}
\end{equation}

From (\ref{triv49}) -- (\ref{triv52}), it follows the differential
equation
\begin{equation}
x\psi_n^{\prime\prime}+(\gamma-2x^2)\psi_n^{\prime}+(2nx-
\frac{\theta_n}{x})\psi_n=0,\qquad{n\geq0},
\label{triv59}
\end{equation}
where $\theta_n=\ds\frac{\gamma(1-(-1)^n)}{2}.$
This equation coincides with the known differential equation
for generalized Hermite  polynomials (\cite{Chih} see also \cite{sze}).
\begin {rem}
For the generators  $a_{\mu}^{+},a_{\mu}^{-}$ of the Heisenberg
algebras $A_{\mu}$, corresponding to the system of the generalized
Hermite polynomials, we have (see \cite{B}):
\begin{equation}
[a_{\mu}^{-},a_{\mu}^{+}]=(\gamma+1)I-2(2B(N)-N). \label{triv60}
\end{equation}
"The energy levels" for Hamiltonian of the related oscillator are
equal to:
\begin{equation}
\lambda_0=\gamma+1,\qquad \lambda_n=2n+\gamma+1,\qquad {n\geq1}.
\label{triv61}
\end{equation}
Finally, we  have for the momentum operator the following formula:
\begin{equation}
P_{\mu}=\imath(\frac{d}{dx}+X_{\mu}^{-1}\Theta_N-X_{\mu}),
\label{triv62}
\end{equation}
\end {rem}
where
\begin{equation}
\Theta_N=2B(N)-N.
\label{trivs61}
\end{equation}

\section {Conditions under which the annihilation operator $A$
can be realized by a differential operator}\label{sA}

Now we consider the relations (\ref{bd2.5}) and (\ref{bd2.6})
without the supposition (\ref{bd2.8})relative to the coefficients
$a_{k, s}$.

Let us introduce
\begin{equation}\label{bd2.21}
C_{k+2m-2w,k+2p+1-2w}= \frac{\left( k+2p+1-2m\right)!}
{b_0^{k+2p+1-2m}}a_{k+2m-2w,k+2p+1-2w}
\end{equation}

Then the following assertion is hold.
\begin{thm}
The operator $A$ defined by the  formula (\ref{bd2.2}) satisfy
the relation (\ref{bd2.3}) iff the elements $a_{k, s}$  of this
operator matrix satisfy the following conditions:
\begin{itemize}
\item
For all  $t=2l,\,n=2p+2l+1,\,p\geq 0,\,l\geq 0$ the following
relation
\begin{equation}
\label{bd2.20}
\begin{array}{c}
\sum_{w=0}^{l-1}\sum_{m=0}^w\left( -1\right) ^m \frac{\alpha
_{2m-1,2p+2l}}{\left( 2w-2m\right) !} B_1 +
\sum_{w=l}^{p+l}\sum_{m=w-l}^w \left(-1\right) ^m \frac{\alpha
_{2m-1,2p+2l}}{\left( 2w-2m\right) !} B_2 =\\ {} \\
=(-1)^p\frac{\left\{ 2l+2p+1\right\} !!}{b_0^{2l}},
\end{array}
\end{equation}
where
\begin{gather*}
B_1=\left(  C_{2l+2m-2w,2p+2l+1-2w}+ \frac{\left(
2p+2l+1-2m\right)}{b_0}
 C_{2l+2m-2w-1,2p+2l-2w}\right) , \\[10pt]
B_2=\left( C_{2l+2m-2w,2p+2l+1-2w}+ \frac{\left(
2p+2l+1-2m\right)}{b_0}
 C_{2l+2m-2w-1,2p+2l-2w}\right)
\end{gather*}
must be hold.
\item
For all $t=2l+1,\,n=2p+2l+2,\,p\geq 0,l\geq 0$
the following relation
\begin{equation}
\label{bd2.22}
\begin{array}{c}
\sum_{w=0}^{l-1}\sum_{m=0}^w\left( -1\right) ^m \frac{\alpha
_{2m-1,2p+2l+1}}{\left( 2w-2m\right) !} B_3 +
\sum_{w=l}^{p+l}\sum_{m=w-l}^w\left( -1\right) ^m B_4 \frac{\alpha
_{2m-1,2p+2l+1}}{\left( 2w-2m\right) !} = \\{}\\
(-1)^p\left\{2l+2p+2\right\} \frac{\alpha
_{2p-1,2p+2l}}{b_0^{2l+1}}
\end{array}
\end{equation}
should be fulfilled, where
\begin{gather*}
B_3=\left( C_{2l+2m-2w+1,2p+2l+2-2w}+ \frac{\left(
2p+2l+2-2m\right) } {\left(2w+1-2m\right)
b_0}C_{2l+2m-2w,2p+2l+1-2w}\right) , \\[10pt] B_4=\left(
C_{2l+2m-2w+1,2p+2l+2-2w}+\frac{\left( 2p+2l+2-2m\right) } {\left(
2w+1-2m\right) b_0}C_{2l+2m-2w,2p+2l+1-2w}\right),
\end{gather*}
\end{itemize}
and the coefficients $\alpha _{2mp-1,n-1}$ are given by (\ref{b13}).
\end{thm}

\section{Conclusion}
From the results given above it follows that in the general case a
orthonormal polynomial system satisfies differential (or
difference) equation of a finite order. Therefore it is
interesting to describe such orthonormal polynomial system, for
which one can obtain a differential (or difference) equation of a
finite order. It is desirable that we have to deal only with
differential equation of the second order. The elementary example
of such class of polynomials selected by the condition (\ref
{bd2.8}), is given by Hermite - Chihara polynomials
\cite{Chih,sze} considered in detail in \cite{Bor}.

Unfortunately, the relations (\ref{bd2.21} -- \ref{bd2.22}) are
rather complicated in the general case. Therefore generally local
formulas representing an annihilation operator (and also for a
momentum operator and Hamiltonian) by a differential operator are
too cumbersome.

Finally, we note the that the formalism put forward in the present
work holds true also after replacement usual derivation by some
of its $q$-analog. The last circumstance provides a useful guide
to investigating of deformed polynomials.

\end{document}